
\baselineskip=14pt
\parskip=10pt

\magnification=\magstephalf
\def\P{{\cal P}}

\def\1{{\overline{1}}}
\def\2{{\overline{2}}}
\parindent=0pt
\overfullrule=0in
\def\Tilde{\char126\relax}
\def\frac#1#2{{#1 \over #2}}
\bf
\centerline
{Automatic Generation of  Generating Functions for Chromatic Polynomials}
\centerline
{
for Grid Graphs (and more general creatures) of Fixed (but arbitrary!) Width
}
\rm
\bigskip
\centerline
{By {\it Shalosh B. EKHAD}, {\it Jocelyn QUAINTANCE}, and {\it Doron ZEILBERGER}}

\bigskip
{\it
\qquad  D\'edi\'e \`a la m\'emoire de 
{\bf PHILIPPE FLAJOLET} (
$1^{er}$ d\'ecembre 1948 - 22 mars 2011)
}
\bigskip
Exclusively published in the {\bf Personal Journal of Shalosh B. Ekhad and Doron Zeilberger}: {\tt http://www.math.rutgers.edu/\Tilde zeilberg/pj.html}
and {\tt arviv.org}.
Written: March 30, 2011.
\bigskip
{\bf Very Important}: This article comments on the Maple package \hfill\break
{\tt http://www.math.rutgers.edu/\Tilde zeilberg/tokhniot/KamaTzviot}. 
Some sample input and output can be gotten from the ``front'' of this article: \hfill\break
{\tt http://www.math.rutgers.edu/\Tilde zeilberg/mamarim/mamarimhtml/tzeva.html} \quad .

\bigskip
One of the {\it many} things done in, or inspired by, the famous INRIA {\bf ALGORITHMS PROJECT} \hfill\break
(http://algo.inria.fr/), the brain child of our beloved guru {\bf Philippe Flajolet},
was to devise efficient algorithms for computing {\it chromatic polynomials} for {\it grid graphs},
namely the Cartesian product $\P_n \times \P_m$ where $\P_n$ is the path of length $n$.
It is fairly easy to see that for a fixed $m$, the generating function
$$
F_m(z,c)=\sum_{n=0}^{\infty} P_{\P_n \times \P_m}(c) z^n 
$$
is a {\it rational function} of {\bf both} $c$ and $z$. 
The Maple package {\tt KamaTzviot} automatically
computes these rational function for {\it any} inputted numeric $m$.
(See procedure {\tt GFk(m,t,c)} of  {\tt KamaTzviot}).

In fact we  do something much more general.
For {\it any} graph $G$,  {\tt KamaTzviot} can (explicitly!) compute

$$
F_G(z,c)=\sum_{n=0}^{\infty} P_{\P_n \times G}(c) z^n \quad .
$$
(See procedure {\tt GFG(G,t,c)} .)

In fact we  do something {\it even} more general!
For {\it any} graph $G$, on $m$ vertices, and for {\it any} bipartite $(m,m)$ graph $C$, let
$M_n(G,C)$ be the graph on $mn$ vertices where  the edges among
$$
{1+im \, , \, 2+im \, , \,  ... \, , \, m+im}
$$
mimic the graph $G$ (for $i=0, \dots , n-1$), and in addition the
edges between \vfill\eject
$$
{1+im \, , \, 2+im \, , \,  ... \, , \,  m+im}
$$
and
$$
{1+(i+1)m \, , \, 2+(i+1)m \, , \,  ... \, , \, m+(i+1)m}
$$
($0 \leq i < n-1$)
mimic the edges of $C$, given as a set of (up to $m^2$) ordered pairs $\{[\alpha,\beta]\}$.
$[\alpha,\beta] \in C$ means that
there is an edge between vertex $\alpha +im$ and vertex $\beta+(i+1)m$ for $0\leq i<n-1$.
Note that when $C$ is the monogamy bipartite 
graph $\{[1,1], \dots, [m,m] \}$,
where Mr $i$ is connected to Mrs $i$ (but no cheating!), then 
$M_n(G,C)$ reduces to the Cartesian product $G \times \P_n$.

{\tt KamaTzviot} can (explicitly!) compute the rational function
(of $z$ and $c$):
$$
F_{G,C}(z,c)=\sum_{n=0} P_{M_n(G,C)}(c) z^n \quad .
$$
See procedure {\tt GFGG(G,C,t,c)}.

{\bf The Method}

Of course we use the {\it transfer matrix method}, but in
a {\it symbolic} context.  The graph $G$ has $m$ vertices, so it
can be legally colored by at most $m$ colors. Let's label the
vertices of $G$, once and for all, by the integers $\{1, \dots , m \}$,
and visualize them from left to right. For each (legal) vertex-coloring of
$G$ we can associate a {\it canonical} form, by renaming the
color of vertex $1$, color $1$, 
then the color of the smallest vertex of $G$ that is not colored by 
color $1$, color $2$, and then the color of the smallest vertex
colored by neither of these colors, color $3$, etc.

For example, If the coloring is $351132$ then its  canonical form is
$123314$.

Our set of ``states''  consists of all possible canonical colorings
of $G$. Of course, since $G$ has finitely many vertices, this
set is a finite set. (In {\tt KamaTzviot} this is done by procedure 
{\tt CC(G)}, {\tt CC} stands for Canonical Colorings). Note that if
$G$ has no edges, then $CC(G)$ is in bijection with the set of set-partitions
of $\{1, \dots, m\}$, so an upper bound for the number of states is the Bell number $B_m$.

Suppose that the ``bottom'' $G$  in $M_n(G,C)$ is in a certain state $S$, and we want to add
another ``layer'' (a copy of $C$ and $G$) to form $M_{n+1}(G,C)$, and we want the
state of the new bottom to be $T$. In how many ways can we color the new $m$ vertices
(namely vertices $mn+1, \dots , mn+m$)
with $c$ colors, so that the new coloring is
still a legal vertex-coloring? If $T$ is the list $[j_1 , j_2, \dots , j_m]$ (of course $j_1=1$),
let's call the {\it actual} colors colored by these vertices $i_{j_1}, \dots, i_{j_m}$ respectively,
and suppose that there are $k$ different colors in that last layer, i.e.
$$
\{ j_1, \dots , j_m \} =\{1, \dots , k \} \quad ,
$$
so that
$$
\{ i_{j_1}, \dots , i_{j_m} \} =\{i_1, \dots , i_k \} \quad .
$$

By the definition of ``state'' any such coloring would not introduce an edge {\it within} this
last layer connecting vertices with the same color, but we have also to worry
about the edges of the last installment of $C$. If the state $S$ is
$$
[a_1, \dots , a_m] \quad ,
$$
and it has $s$ different colors, i,e,
$$
\{ a_1, \dots , a_m \} =\{1, \dots , s \} \quad ,
$$
then for any edge $[\alpha, \beta] \in C$, we know, by the construction of $M_n(G,C)$,
that there is an edge between $mn+\alpha$ and $m(n+1)+ \beta$.
This entails that $i_{j_\beta}$ can {\bf not}
be $a_{\alpha}$. So for each and every $i_1, \dots, i_k$ there is a set of ``forbidden colors'',
out of the $s$ colors of layer $n$. Of course each of $i_1, \dots, i_k$
can also be colored in one of the $c-s$ colors that are not used in layer $n$.
Converting the ``negative'' conditions into ``positive'' ones, we get the set of permissible
colors for each $i_1, \dots, i_k$,  including what we called ``option $0$'' (in procedure {\tt TS1S2G} of our Maple package)
that denotes choosing one of the remaining $c-s$ colors. Taking the Cartesian product
of the option-sets for each of $i_1, \dots, i_k$,  we get atomic events where
some of the members of $\{ i_1, \dots , i_k \}$
are committed to be one of the $s$ colors of layer $n$, and the rest are
{\it different} colors from those $c-s$ colors. If there are $\gamma$ such ``$0$''s, of course
the number of ways of doing it is the polynomial of degree $\gamma$ in $c$, $\gamma!{{c-s \choose \gamma}}$.
Adding the number of possibilities of all the atomic options would give the matrix-entry connecting state $S$ to state $T$.

We (or rather the first-named author) does it for each pair of states $S$ and $T$, building up the
{\it transfer matrix} {\bf completely automatically} (in other words, it does the ``combinatorial research''
all by itself!). Once we have the transfer matrix, Ekhad sets up the obvious set of equations
for the generating functions for colorings ending  at any given state, solves, (symbolically and automatically!)
the resulting set of linear equations (with coefficients that are poynomials of $z$ and $c$), and then adds them up to get the desired generating function.

Let us conclude with a simple example of the grid graphs $P_3 \times P_n$. Here
$$
Edges(G)=\{ \{1,2\}, \{2, 3 \} \} \quad ,
$$
$$
C=\{ [1,1],[2,2],[3,3] \} \quad .
$$
There are two states: $121$ and $123$. Let's compute the matrix entry
connecting state $123$ to $121$.

The coloring of the bottom layer is $[i_1 ,i_2 ,i_1]$ for {\it different} colors $i_1$ and $i_2$ (chosen between $1$ and $c$).

Because of the edge $[1,1]$ of $C$ we have the restriction $i_1 \neq 1$.

Because of the edge $[2,2]$, we have the restriction $i_2 \neq 2$.

Because of the edge $[3,3]$, we have the restriction $i_1 \neq 3$.

So $i_1$ may {\bf not} be colored with color $1$ and color $3$, while $i_2$ may
{\bf not} be colored with color $2$. Going to the {\bf positive} rephrasing:

$i_1=2$ OR $3 \leq i_1 \leq c$

AND

$i_2=1$ OR $i_2=3$ OR $3 \leq i_2 \leq c$ \quad .

Taking the ``product'' we have the six ``atomic'' events:

$i_1=2$ AND $i_2=1$ ($1$ possibility)

OR

$i_1=2$ AND $i_2=3$ ($1$ possibility)

OR

$i_1=2$ AND $3 \leq i_2 \leq c \,\,$  ($1!{{c-3} \choose {1}}=c-3$ possibilities)

OR

$3 \leq i_1 \leq c$ AND $i_2=1 \,\, $ ($1!{{c-3} \choose {1}}=c-3$ possibilities)

OR

$3 \leq i_1 \leq c$ AND $i_2=3\,\,$ ($1!{{c-3} \choose {1}}=c-3$ possibilities)

OR

$3 \leq i_1 \leq c$ AND $3 \leq i_2 \leq c \,\,$  
($2!{{c-3} \choose {2}}=(c-3)(c-4)$ possibilities) .

Adding these up gives the matrix entry:
$$
M[123,121]=2 \cdot 1 + 3(c-3)+ (c-3)(c-4)=c^2-4c+5 \quad .
$$

We leave to our {\it human} readers, as an instructive exercise to test their comprehension of our method,
to verify that
$$
M[121,121]=c^2-3c+3 \quad , \quad M[121,123]=c^3-6c^2+13c-10 \quad , \quad M[123,123]=c^3-6c^2+14c-13 \quad .
$$
Of course, {\tt KamaTzviot} can do so much {\bf more}.

\quad

\hrule
Shalosh B. Ekhad, Jocelyn Quaintance, and Doron Zeilberger,
Department of Mathematics, Hill Center-Busch Campus, Rugters University, 110 Frelinghuysen Rd,
Piscataway, NJ 08544, USA.

Email: {\tt [c/o zeilberg, quaintan, zeilberg]@math.rutgers.edu} .

\end